\documentclass[12pt]{amsart}
\usepackage[margin=2.0cm]{geometry}

\usepackage{amssymb}
\usepackage{amsfonts}
\usepackage{mathrsfs}
\usepackage{cite}
\usepackage{graphicx}

\usepackage{pgfplots}
\pgfplotsset{compat=1.15}
\usepackage{mathrsfs}
\usetikzlibrary{arrows}
\usepackage{float}

\definecolor{ffqqqq}{rgb}{1,0,0}
\definecolor{qqzzqq}{rgb}{0,0.6,0}

\newcommand{\R}{{\mathbb R}}
\newcommand{\N}{{\mathbb N}}
\newcommand{\hyp}{{\mathbb H}}

\newtheorem{theorem}{Theorem}[section]

\newtheorem{lemma}[theorem]{Lemma}
\newtheorem{prop}[theorem]{Proposition}
\newtheorem{corollary}[theorem]{Corollary}

\newcommand{\Int}[1]{\text{\rm int }#1}

\newcommand{\bigslant}[2]{{\raisebox{.2em}{$#1$}\left/\raisebox{-.2em}{$#2$}\right.}}

\DeclareMathOperator{\diam}{\mathrm{diam}}

\DeclareMathOperator{\conv}{\mathrm{conv}}

\makeatletter
\DeclareFontFamily{U}{tipa}{}
\DeclareFontShape{U}{tipa}{m}{n}{<->tipa10}{}
\newcommand{\arc@char}{{\usefont{U}{tipa}{m}{n}\symbol{62}}}%

\newcommand{\arc}[1]{\mathpalette\arc@arc{#1}}

\newcommand{\arc@arc}[2]{%
  \sbox0{$\m@th#1#2$}%
  \vbox{
    \hbox{\resizebox{\wd0}{\height}{\arc@char}}
    \nointerlineskip
    \box0
  }%
}
\makeatother

\title[Bodies of constant width in the hyperbolic space]{Hyperbolic width functions and characterizations of bodies of constant width in the hyperbolic space}
\author{K\'aroly J. B\"or\"oczky, Andr\'as Cs\'epai, \'Ad\'am Sagmeister }
\address{Alfr\'ed R\'enyi Institute of Mathematics,  Realtanoda u. 13-15, H-1053 Budapest, Hungary} \email{boroczky.karoly.j@renyi.hu}
\address{E\"otv\"os Lor\'and University, Institute of Mathematics, P\'azm\'any P\'eter s\'et\'any 1/c, Budapest, H-1117 Hungary} \email{csepai.andras112358@gmail.com }
\address{E\"otv\"os Lor\'and University, Institute of Mathematics, P\'azm\'any P\'eter s\'et\'any 1/c, Budapest, H-1117 Hungary} \email{sagmeister.adam@gmail.com }
\subjclass[2010]{Primary: }

\begin{document}

\maketitle

\begin{abstract}
We discuss 
basic properties of several different 
width functions in the $n$-dimensional hyperbolic space such as continuity, and we also define a new hyperbolic width as the extension of Leichtweiss' width function. Then we prove a characterization theorem of bodies of constant width regarding the aforementioned notions of hyperbolic width.
\end{abstract}

\section{Introduction}

Let $\hyp^n$ be the $n$-dimensional hyperbolic space for $n\geq 2$. We write $d(x,y)$ to denote the geodesic distance between $x,y\in \hyp^n$ and $[x,y]$ to denote the geodesic segment between $x$ and $y$ whose length is $d(x,y)$. For $X\subseteq \hyp^n$, we say that $X$ is \emph{convex} if $[x,y]\subseteq X$ for $x,y\in X$. In addition, $X$ is a \emph{convex body}, if $X$ is compact and convex with non-empty interior. For $z\in \hyp^n$ and $r\geq 0$, we write
$B(z,r)=\{x\in \hyp^n\colon d(x,z)\leq r\}$ to denote the ball of radius $r$ centered at $z$.

Bodies of constant width play an important role in the well-researched field of convex geometric inequalities in the Euclidean space or in Minkowski spaces (see Martini, Swanepoel\cite{MS04}). Some more recent studies investigate bodies of constant width on the sphere (see Lassak, Musielak\cite{LaM18}). However, when it comes to the hyperbolic space, there are many different attempts to define the width function of a convex body (see Santal\'o\cite{S45}, Fillmore\cite{Fil70}, Leichtweiss\cite{Lei05}, Jer\'onimo-Castro--Jimenez-Lopez\cite{JCJL17}, G. Horv\'ath\cite{Hor21}). We will also introduce a new concept of hyperbolic width, based on Leichtweiss' definition. Sometimes we cannot even compare these functions, and even when we can, they often take different values. The goal of this paper is to show, that this difference disappears within the class of bodies of constant width.

\medskip\noindent\textbf{Organisation of the paper.}
In Section \ref{models} we introduce the models of the hyperbolic space we will use throughout the paper. In Section \ref{cptsubsets} we describe basic properties of compact subsets in $\hyp^n$. In Section \ref{width} we define the width functions $w_S$, $w_F$, $w_L$, $w_{JCJL}$, $w_{GH}$ of Santal\'o, Fillmore, Leichtweiss, Jer\'onimo-Castro--Jimenez-Lopez and G. Horv\'ath, respectively, and introduce the new width function $w$. Section \ref{cont} is about the continuity of these width functions, in particular we prove:

\begin{theorem}
\label{contthm}
The width functions $w_S(K,(z,H_z)$, $w_F(K,i)$, $w_L(K,(p,H_p))$, $w_{JCJL}(K,(z,\ell_z,z'))$, $w_{GH}(K,i)$ and $w(K,H)$ are continuous in all of their parameters. The minimal widths $w_S(K)$, $w_F(K)$, $w_L(K,p)$, $w_{JCJL}(K)$, $w_{GH}(K)$ and $w(K)$ of convex bodies $K$ are also continuous in $K$.
\end{theorem}


In Section \ref{diam} we verify that the diameter coincides with the maximal width using any of the width notions defined in the paper.


Section \ref{constwidth} deals with complete bodies and bodies of constant width. It contains our main result:

\begin{theorem}
\label{completeness}
Let
$K\subset \hyp^n$ be a convex body of diameter $D>0$. The following are equivalent:
\begin{description}
\item[(i)] $K$ is complete;
\item[(ii)] $K=\bigcap_{x\in K}B\left(x,D\right)$;
\item[(iii)] $K$ is of constant width $D$;
\item[(iv)] $w_S\left(K,\left(z,H_z\right)\right)=D$ for all $z\in\partial K$ and $H_z$ supporting hyperplanes to $K$ at $z$;
\item[(v)] $w_F\left(K,i\right)=D$ for all ideal points $i$;
\item[(vi)] $w_L\left(K,\left(p,H_p\right)\right)=D$ for all hyperplanes $H_p$ through a fixed inner point $p\in\Int K$;
\item[(vii)] $w\left(K,H\right)=D$ for all hyperplanes $H$ intersecting $K$.
\end{description}
\end{theorem}

Finally, in Section \ref{balls} we consider the circumscribed and inscribed balls of convex bodies and deduce results on the circumradii and inradii of complete bodies. In particular, we prove a hyperbolic version of Scott's theorem:

\begin{theorem}\label{scott}
Let $K\subset \hyp^n$ be a convex body of diameter $D$ and circumscribed ball $B$. Then, there exists a completion $\widetilde{K}$ of $K$ such that $B$ is also the circumscribed ball of $\widetilde{K}$.
\end{theorem}

\section{Models of the hyperbolic space}
\label{models}

\subsection{The hyperboloid model}

We assume that $\hyp^n$ is embedded into $\R^{n+1}$ using the hyperboloid model. We write $\langle \cdot,\cdot\rangle$ to denote the standard scalar product in $\R^{n+1}$, and write 
$z^\bot=\{x\in\R^{n+1}\colon\langle x,z\rangle=0\}$ for a $z\in\R^{n+1}\setminus\left\{o\right\}$. Fix an $e\in S^n$, then we have
$$
\hyp^n=\left\{x+te\colon x\in e^\bot \mbox{ and }t\geq 1\mbox{ and }t^2-\langle x,x\rangle=1\right\}.
$$
We also consider the following symmetric bilinear form $\mathcal{B}$ on $\R^{n+1}$: If $x=x_0+te\in\R^{n+1}$ and $y=y_0+se\in\R^{n+1}$ for
$x_0,y_0\in e^\bot$ and $t,s\in\R$, then
$$
\mathcal{B}(x,y)=ts-\langle x_0,y_0\rangle.
$$
In particular,
\begin{equation}
\label{Hn}
\mathcal{B}(x,x)=1\mbox{ for $x\in \hyp^n$}.
\end{equation}
For $z\in \hyp^n$, we define the tangent space $T_z$ as
$$
T_z=\left\{x\in\R^{n+1}\colon\mathcal{B}(x,z)=0\right\}.
$$
We observe that $T_z$ is an $n$-dimensional real vector space equipped with the scalar product $-\mathcal{B}(\cdot,\cdot)$.

For $z\in \hyp^n$ and unit vector $u\in T_z$, the geodesic line $\ell$ passing through $z$ and determined by $u$ consists of the points
$$
p_t=z\,{\rm cosh}\,t+u\,{\rm sinh}\,t
$$
for $t\in \R$. Here the map $t\mapsto p_t$ is bijective onto $\ell$ and satisfies $d(z,p_t)=|t|$ for $t\in \R$. If $t>0$, then we say that $u$ points towards $p_t$ along the geodesic segment
$$
[z,p_t]=\{p_s\colon 0\leq s\leq t\}
$$
of length $t$.

A hyperplane $H$ in $\hyp^n$ passing through the point $z\in \hyp^n$ and having exterior unit normal $u\in T_z$, and the corresponding half-spaces $H^-$ and $H^+$ are defined as follows: $H^-=\hyp^n\setminus{\rm int} H^+$, $H=\{x\in \hyp^n\colon\mathcal{B}(x,u)=0\}$ and $H^+=\{x\in \hyp^n\colon-\mathcal{B}(x,u)\geq 0\}$.

In $\hyp^n$, for a hyperplane $H$, the \emph{hyperball} of radius $\varepsilon\geq 0$ is the \emph{parallel domain} $H^{\left(\varepsilon\right)}$, that is the set of points in $\hyp^n$ of distance less than or equal to $H$. The boundary components of a hyperball are called \emph{hyperspheres} (or \emph{hypercycles} in the 2-dimensional case). A hypersphere $\partial\mathcal{H}$ of the hyperball $H$ has two connected components, $H_1$ and $H_2$ (except the degenerate case where $\varepsilon=0$ and $H=H_1=H_2$). Both $H_1$ and $H_2$ divide the hyperbolic space $\hyp^n$ into two closed connected components, $H_k^+$ and $H_k^-$ for $k\in\left\{1,2\right\}$, such that $H_k^+\cap H_k^-=H_k$, $H_1^+\cap H_2^+$ is convex, and $H_k^-\subseteq H_{3-k}^+$.

We also consider horoballs and horospheres ``centered'' at ideal points. An \emph{ideal point} of $\hyp^n$ is represented by the linear hull of a $z\in\R^{n+1}$ with $\mathcal{B}(z,z)=0$ and $\mathcal{B}(z,e)>0$. Any hyperbolic line $\ell=\Pi\cap \hyp^n$ for a linear two-plane $\Pi$ with $\Pi\cap \hyp^n\neq\emptyset$ contains exactly two ideal points. If $p\in\ell$ and $v\in T_p\cap \Pi$ is one of the two tangent vectors to $\ell$ at $p$ with $\mathcal{B}(v,v)=-1$, then the two ideal points are represented by the spans of $u-v$ (``point of $\ell$ at infinity'' in the direction of $v$) and $u+v$ (``point of $\ell$ at infinity'' in the direction of $-v$). We also observe that hyperbolic lines containing either ideal point of $\ell$ are the lines parallel to $\ell$.

Let $z\in \R^{n+1}$ with $\mathcal{B}(z,z)=0$ and  
$\mathcal{B}(z,e)>0$, and let $i$ be the ideal point represented by ${\rm lin}\,z$. For any $p\in \hyp^n$, there exists a unique line $\ell$ passing through $p$ and $i$, and a unit tangent vector $v\in T_p$ to $\ell$ is 
$u-\mathcal{B}(p,z)^{-1}z$.

For $s>0$, a \emph{horoball} at $i$ is
$$
A=\{x\in \hyp^n\colon\mathcal{B}(z,x)\geq s\},
$$
and the corresponding \emph{horosphere} is
$$
\partial A=\{x\in \hyp^n\colon\mathcal{B}(z,x)=s\}.
$$
We observe that for any $p\in\partial A$, the line passing through $p$ and $i$ is orthogonal to the horosphere $\partial A$. If $n=2$, then a horoball is also called a \emph{horodisk} and a horosphere is called a \emph{horocycle} as well.

\subsection{The Poincar\'e ball model}

We can obtain the Poincar\'e ball model by a radial projection of the hyperboloid through the apex of the hyperboloid $e\in\R^{n+1}$ onto $e^{\perp}$. In the Poincar\'e ball model, the hyperbolic space $\hyp^n$ is identified with the interior of the unit Euclidean ball $B^n$ in $\R^n$, and the set of ideal points are just $\partial B^n$. 
A hyperbolic line in the Poincar\'e ball model is the intersection of ${\rm int}B^n$ and a Euclidean circle that is orthogonal to $\partial B^n$ at the two intersection points.
A hyperplane is the intersection of a Euclidean $(n-1)$-sphere with ${\rm int}B^n$ orthogonal to $\partial B^n$ at the $(n-2)$-sphere of intersection points.
A horosphere at an ideal point $i\in \partial B^n$
is of the form $\partial G\setminus\left\{i\right\}$ for a Euclidean $n$-ball $G\subset B^n$ of radius less than one and touching $\partial B^n$ in $i$.
In addition, hyperbolic $n$-balls, in the Poincar\'e model coincide with Euclidean $n$-balls contained in ${\rm int}\,B^n$.

\subsection{The Beltrami--Cayley--Klein model}

The Beltrami--Cayley--Klein model can be derived from the hyperboloid model by a radial projection onto $e^{\perp}-e$ through the origin $o\in\R^{n+1}$. Here, the points of the hyperboloid are the interior points of the Euclidean unit ball centered at $-e$, while ideal points are again the points of this unit ball. Sometimes it is easier to investigate convexity, as hyperplanes and lines of $\hyp^n$ are intersections of the ball with Euclidean hyperplanes and lines. In particular, hyperbolic convexity coincides with Euclidean convexity in the Beltrami--Cayley--Klein ball.

\section{Compact subsets in the hyperbolic space}
\label{cptsubsets}

We impose a metric on compact subsets. For a compact set $C\subset \hyp^n$ and $z\in \hyp^n$, we set $d(z,C)=\min_{x\in C}d(z,x)$. For any non-empty compact set $C_1,C_2\subset \hyp^n$, we define their \emph{Hausdorff distance}
$$
\delta(C_1,C_2)=\max\left\{\max_{x\in C_2}d(x,C_1),\max_{y\in C_1}d(y,C_2)\right\}.
$$
The Hausdorff distance is a metric on the space of compact subsets in $\hyp^n$. We say that a sequence $\left\{C_m\right\}_{m\in\N}$ of compact subsets of $\hyp^n$ is \emph{bounded} if there is a ball containing every $C_m$.
For compact sets $C_m,C\subset \hyp^n$, we write $C_m\to C$ to denote if the sequence $\left\{C_m\right\}_{m\in\N}$ tends to $C$ in terms of the Hausdorff distance.

The following Lemma~\ref{ballinhoroball} is  a simple consequence of these properties of the Poincar\'e disk model.

\begin{lemma}
\label{ballinhoroball}
Any horoball $A$ in $\hyp^n$ is closed and convex, and if $A$ and a ball $B(z,r)$, $r>0$, have a common exterior normal $u\in T_y$ for some $y\in \partial A\cap \partial B(z,r)$, then
$B(z,r)\subset A$ with $B(z,r)\cap A=\{y\}$.
\end{lemma}

Lemma~\ref{ballinhoroball} directly yields the following.

\begin{corollary}\label{horostripwidth}
If $A$ and $C$ are horoballs at an ideal point $i$ of $\hyp^n$ with $C\subset A$, then  $d(x,y)\geq\Delta\left(A,C\right)$ for any $x\in\partial A$ and $y\in\partial C$ with equality if and only if the line $\ell$ passing through $x$ and $y$ has $i$ as an ideal point (and hence $\ell$ is orthogonal to  $\partial A$ and $\partial C$). Here $\Delta\left(A,C\right)$ denotes the lenth of the geodesic segment obtained by the intersection of the line through the origin and $i$, and the horospherical domain $A\setminus\Int C$.
\end{corollary}

The following is well-known  (see B\"or\"oczky, Sagmeister \cite{BoS20}).

\begin{lemma}
\label{Hausdorffconvergence}
For compact sets $C_m,C\subset \hyp^n$, we have $C_m\to C$ if and only if
\begin{description}
\item[(i)] assuming $x_m\in C_m$, the sequence $\{x_m\}$ is bounded and any accumulation point of  $\{x_m\}$  
lies in $C$;
\item[(ii)] for any $y\in C$, there exist $x_m\in C_m$ for each $m$ such that $\lim_{m\to \infty}x_m=y$.
\end{description}
\end{lemma}
The space of compact subsets of $\hyp^n$ is locally compact according to the Blaschke Selection Theorem (see R. Schneider \cite{Sch14}).

\begin{theorem}[Blaschke]
\label{BlaschkeSel}
Any bounded sequence of compact subsets of $\hyp^n$ has a convergent subsequence.
\end{theorem}

For convergent sequences of compact subsets of $\hyp^n$, we have the following 
(see B\"or\"oczky, Sagmeister \cite{BoS20}).

\begin{lemma}
\label{limit}
Let the sequence $\{C_m\}$ of compact subsets of $\hyp^n$ tend to $C$. Then,
$${\rm diam}\,C=\lim_{m\to \infty}{\rm diam}\,C_m.$$
\end{lemma}

Recall that for any set $X\subset \hyp^n$ and $\varrho\geq 0$, the
parallel domain is
$$
X^{(\varrho)}=\{z\in \hyp^n\colon\exists x\in X\mbox{ \ with \ }
d(x,z)\leq \varrho\}=\bigcup \{B(x,\varrho)\colon x\in X\}.
$$
The triangle inequality and considering $x,y\in X$ with 
$d(x,y)={\rm diam}\, X$ show that
\begin{equation}
\label{parallel-diam}
{\rm diam}\, X^{(\varrho)}=2\varrho+{\rm diam}\, X.
\end{equation}
The following is well-known (see B\"or\"oczky, Sagmeister\cite{BoS23}):
\begin{lemma}
\label{limit-parallel}
If $\varrho\geq 0$, and the sequence $\{C_m\}$ of compact subsets of $\hyp^n$ tends to $C$, then
\begin{description}
\item[(i)] $\left\{C_m^{(\varrho)}\right\}$ tends to $C^{(\varrho)}$;
\item[(ii)] ${\rm diam}\,C^{(\varrho)}=\lim_{m\to \infty}{\rm diam}\,C_m^{(\varrho)}$.
\end{description}
\end{lemma}

The following lemma shows that it is true, that any parallel domain of a convex body is also a convex body in the hyperbolic space.

\begin{lemma}\label{parallel_convexity}
Let $K\subset \hyp^n$ be a convex body of diameter $D>0$ and $\varepsilon>0$. Then the parallel domain $K^{\left(\varepsilon\right)}$ is also a convex body.
\end{lemma}

\proof

It is clear that to verify that $K^{\left(\varepsilon\right)}$ is a convex body, we only need to prove its convexity. We use the Beltrami--Cayley--Klein model to verify that. In this model, hyperbolic lines are Euclidean lines, therefore it is enough to check Euclidean convexity in this model. Let $x,y\in K^{\left(\varepsilon\right)}$ and $\lambda\in\left[0,1\right]$. We want to see that $\lambda x+\left(1-\lambda\right)y\in K^{\left(\varepsilon\right)}$. By the definition of the parallel domain, there are points $x_0,y_0\in K$, such that $d\left(x,x_0\right)\leq\varepsilon$ and $d\left(y,y_0\right)\leq\varepsilon$. Let $\ell$ be the line connecting $x$ and $y$. It is well-known, that hypercycles are convex in the Beltrami--Cayley--Klein disk, and $x_0,y_0$ are points in the hypercycle $\ell^{\left(\varepsilon\right)}$, so the geodesic segment $\left[x_0,y_0\right]\subset\ell^{\left(\varepsilon\right)}$. Now let $z\in\left[x_0,y_0\right]$ be an arbitrary point, and let $x_1$ and $y_1$ be the intersection points of the line through $x_0$ and $y_0$, and the orthogonal hyperplane of $\ell$ through the points $x$ and $y$, respectively. We may assume that $\ell$ is a diameter of the disk, so the orthogonal hyperplane is orthogonal in the Euclidean setup as well. If $z\in\left[x_0,y_0\right]\cap\left[x_1,y_1\right]\subset\ell^{\left(\varepsilon\right)}$, then it is clear, that $d\left(z,z'\right)\leq\varepsilon$, where $z$ is the intersection point of $\ell$ and the orthogonal hyperplane to $\ell$ through $z$. The only questionable case is when $z\in\left[x_0,x_1\right]$ or $z\in\left[y_0,y_1\right]$. We may assume that the first is true. But then $z$ is on the side $\left[x_0,x_1\right]$ of the hyperbolic triangle $\left[x,x_0,x_1\right]$, and both of the sides $\left[x,x_0\right]$ and $\left[x,x_1\right]$ are at most $\varepsilon$ long, and hence $d\left(x,z\right)\leq\varepsilon$ in this case.

\endproof

\section{Hyperbolic width functions}
\label{width}

\subsection{Santal\'o's width function}

The following width function was introduced by Santal\'o in $\hyp^2$, but has a natural extension on higher dimensions. Let $K\subset \hyp^n$ be a convex body. Let $H_z$ be a supporting hyperplane to $K$ at a boundary point $z\in\partial K$ and $\ell$ be the unique hyperbolic line orthogonal to $H_z$ through $z$. Then there is a unique supporting hyperplane $H'_z$ to $K$, which is orthogonal to $\ell$ and is disjoint from $H_z$. If $\left\{z'\right\}=H'_z\cap\ell$, we call $d\left(z,z'\right)$ the \emph{Santal\'o width} of $K$ at $\left(z,H_z\right)$, and we use the notation $$w_S\left(K,\left(z,H_z\right)\right).$$

The following example shows that although Santal\'o's definition is a very natural extension of the Euclidean width function, it is not monotonous, considering the minimal width. Let
$$
w_S\left(K\right)=\min\left\{w_S\left(K,\left(z,H_z\right)\right)\colon z\in\partial K\text{ and }H_z\text{ is a supporting hyperplane of }K\text{ at }z\right\}.
$$
\begin{center}
\begin{figure}[h]
\begin{tikzpicture}[line cap=round,line join=round,>=triangle 45,x=1cm,y=1cm,scale=3.5]
\clip(-1.1,-1.1) rectangle (1.1,1.1);
\draw [line width=1pt,dashed] (0,0) circle (1cm);
\draw [shift={(-0.8660254037844384,0.5)},line width=2pt]  plot[domain=-1.0471975511965983:0,variable=\t]({1*0.8660254037844385*cos(\t r)+0*0.8660254037844385*sin(\t r)},{0*0.8660254037844385*cos(\t r)+1*0.8660254037844385*sin(\t r)});
\draw [shift={(0,-1)},line width=2pt]  plot[domain=1.0471975511965974:2.0943951023931957,variable=\t]({1*0.8660254037844383*cos(\t r)+0*0.8660254037844383*sin(\t r)},{0*0.8660254037844383*cos(\t r)+1*0.8660254037844383*sin(\t r)});
\draw [shift={(0.866025403784439,0.5)},line width=2pt]  plot[domain=3.141592653589793:4.1887902047863905,variable=\t]({1*0.866025403784439*cos(\t r)+0*0.866025403784439*sin(\t r)},{0*0.866025403784439*cos(\t r)+1*0.866025403784439*sin(\t r)});
\draw [shift={(0,-1)},line width=1pt,dotted]  plot[domain=0.4478323969289319:1.0471975511965974,variable=\t]({1*0.8660254037844385*cos(\t r)+0*0.8660254037844385*sin(\t r)},{0*0.8660254037844385*cos(\t r)+1*0.8660254037844385*sin(\t r)});
\draw [shift={(0,-1)},line width=1pt,dotted]  plot[domain=2.0943951023931957:2.6937602566608607,variable=\t]({1*0.8660254037844385*cos(\t r)+0*0.8660254037844385*sin(\t r)},{0*0.8660254037844385*cos(\t r)+1*0.8660254037844385*sin(\t r)});
\draw [shift={(0.8396429153503638,1.25)},line width=1pt,dotted]  plot[domain=3.3946252068732403:4.847174744134745,variable=\t]({1*1.1258331249781461*cos(\t r)+0*1.1258331249781461*sin(\t r)},{0*1.1258331249781461*cos(\t r)+1*1.1258331249781461*sin(\t r)});
\draw [shift={(2.1214375390242237,-0.625)},line width=1pt,dotted]  plot[domain=2.3858969368393477:2.689071965658022,variable=\t]({1*1.9725927689163716*cos(\t r)+0*1.9725927689163716*sin(\t r)},{0*1.9725927689163716*cos(\t r)+1*1.9725927689163716*sin(\t r)});
\draw [shift={(2.1214375390242237,-0.625)},line width=2pt,color=red]  plot[domain=2.689071965658022:2.9028591627707554,variable=\t]({1*1.9725927689163711*cos(\t r)+0*1.9725927689163711*sin(\t r)},{0*1.9725927689163711*cos(\t r)+1*1.9725927689163711*sin(\t r)});
\draw [shift={(2.1214375390242237,-0.625)},line width=1pt,dotted]  plot[domain=2.9028591627707554:3.324276450977196,variable=\t]({1*1.9725927689163714*cos(\t r)+0*1.9725927689163714*sin(\t r)},{0*1.9725927689163714*cos(\t r)+1*1.9725927689163714*sin(\t r)});
\draw (-0.11905516326626699,0.17111576568823916) node[anchor=north west] {$K$};
\draw (0.06,-0.14) node[anchor=north west] {$z$};
\draw (0.29,0.2) node[anchor=north west] {$\ell$};
\draw (0.38,0.4) node[anchor=north west] {$z'$};
\draw (-0.5727838405248699,-0.27) node[anchor=north west] {$H_z$};
\draw (-0.08,0.72) node[anchor=north west] {$H_z'$};
\begin{scriptsize}
\draw [fill=black] (0.20479095643453793,-0.15853659368774287) circle (0.6pt);
\draw [fill=black] (0.34739051187511893,0.2374844215663137) circle (0.6pt);
\end{scriptsize}
\end{tikzpicture}
\caption{The Santal\'o width $w_S\left(K,\left(z,H_z\right)\right)$ of a regular triangle $K$ at a boundary point $z$ and supporting line $H_z$ at $z$}
\end{figure}
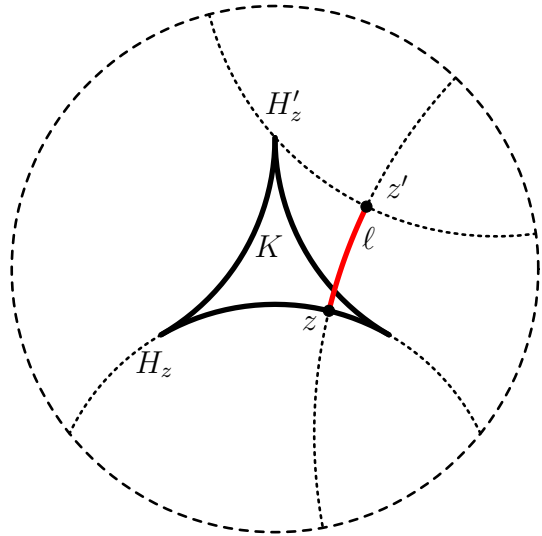
\end{center}

\begin{prop}
There is a convex body $K\subset \hyp^2$ containing the disk $B=B\left(z,r\right)\subset \hyp^2$, such that $w_S\left(K\right)<w_S\left(B\right)$.
\end{prop}
\proof
There are tangent ultraparallel lines $\ell_1$ and $\ell_2$ to the disk $B$, such that the unique perpendicular line $\ell$ to $\ell_1$ and $\ell_2$ both is exterior to $B$. Let $y_1$ be the intersection point of $\ell$ and $\ell_1$, while $y_2$ denotes the intersection point of $\ell$ and $\ell_2$. If $\ell_1$ touches $B$ at $z_1$ and $\ell_2$ touches $B$ at $z_2$, then
$$
w_S\left(B\right)=2r=d\left(z,z_1\right)+d\left(z,z_2\right)>d\left(z_1,z_2\right)>d\left(y_1,y_2\right)=w_S\left(\widetilde{K},\left(y_1,\ell_1\right)\right)\geq w_S\left(\widetilde{K}\right),
$$
where $\widetilde{K}=\conv B\cup\left[y_1,y_2\right]$.
\endproof

\begin{center}
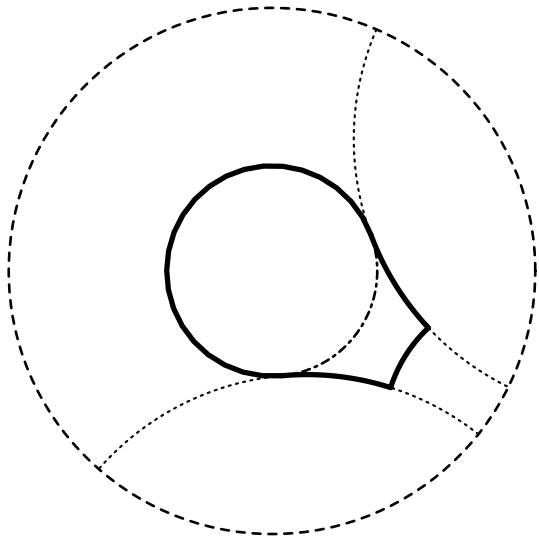
\begin{figure}[h]
\begin{tikzpicture}[line cap=round,line join=round,>=triangle 45,x=1cm,y=1cm,scale=3.5]
\clip(-1.1,-1.1) rectangle (1.1,1.1);
\draw [line width=1pt,dashed] (0,0) circle (1cm);
\draw [shift={(0,0)},line width=2pt]  plot[domain=0.3531920212392142:4.8037816138247535,variable=\t]({1*0.4*cos(\t r)+0*0.4*sin(\t r)},{0*0.4*cos(\t r)+1*0.4*sin(\t r)});
\draw [shift={(0,0)},line width=1pt,dashdotted]  plot[domain=-1.4794036933548327:0.3531920212392142,variable=\t]({1*0.4*cos(\t r)+0*0.4*sin(\t r)},{0*0.4*cos(\t r)+1*0.4*sin(\t r)});
\draw [shift={(0.9579668620340719,-0.6047497503326486)},line width=2pt]  plot[domain=2.3233796344600157:2.8335940006039526,variable=\t]({1*0.5323746512399019*cos(\t r)+0*0.5323746512399019*sin(\t r)},{0*0.5323746512399019*cos(\t r)+1*0.5323746512399019*sin(\t r)});
\draw [shift={(0.13233491507719497,-1.443948569115781)},line width=.8pt,dotted]  plot[domain=0.9011762060102306:1.2627976738090552,variable=\t]({1*1.05*cos(\t r)+0*1.05*sin(\t r)},{0*1.05*cos(\t r)+1*1.05*sin(\t r)});
\draw [shift={(0.13233491507719497,-1.443948569115781)},line width=2pt]  plot[domain=1.2627976738090552:1.6621889602349604,variable=\t]({1*1.05*cos(\t r)+0*1.05*sin(\t r)},{0*1.05*cos(\t r)+1*1.05*sin(\t r)});
\draw [shift={(0.13233491507719497,-1.443948569115781)},line width=.8pt,dotted]  plot[domain=1.6621889602349604:2.4232017144596902,variable=\t]({1*1.05*cos(\t r)+0*1.05*sin(\t r)},{0*1.05*cos(\t r)+1*1.05*sin(\t r)});
\draw [shift={(1.3604964183880752,0.5015471020285323)},line width=.8pt,dotted]  plot[domain=2.7337719206042776:3.4947846748290075,variable=\t]({1*1.05*cos(\t r)+0*1.05*sin(\t r)},{0*1.05*cos(\t r)+1*1.05*sin(\t r)});
\draw [shift={(1.3604964183880752,0.5015471020285323)},line width=2pt]  plot[domain=3.4947846748290075:3.8941759612549123,variable=\t]({1*1.05*cos(\t r)+0*1.05*sin(\t r)},{0*1.05*cos(\t r)+1*1.05*sin(\t r)});
\draw [shift={(1.3604964183880752,0.5015471020285323)},line width=.8pt,dotted]  plot[domain=3.8941759612549123:4.255797429053737,variable=\t]({1*1.05*cos(\t r)+0*1.05*sin(\t r)},{0*1.05*cos(\t r)+1*1.05*sin(\t r)});
\end{tikzpicture}
\caption{The minimal Santal\'o width is not monotonic}
\end{figure}
\end{center}

\noindent {\bf Remark. }Similarly we can construct a counterexample in $\hyp^n$ for any dimension $n\geq 2$.

\subsection{Fillmore's width function}

We call a set \emph{h-convex}, if it is the intersection of some horoballs. A set in $\hyp^n$ is a \emph{h-convex body} if it is h-convex, compact, and of non-empty interior. As horoballs are convex, every h-convex body is also a convex body, but the reverse is not true in general, polytopes for example are convex but not h-convex.

If $r>s$, then the horoball $C=\{x\in \hyp^n\colon\mathcal{B}(z,x)\geq r\}\subset A$ satisfies that the distance of the parallel horospheres $\partial A$ and $\partial C$ is $d>0$ where
$d=\ln\left(\frac{r}{s}\right)$. In particular,  any hyperbolic line $\ell$ containing $i$ satisfies that the distance between 
$\ell\cap \partial A$ and $\ell\cap \partial C$ is $d$, and we set $d=w_{\mathrm{hor}}\left(\partial A,\partial C\right)$
be the width of closed region bounded by the horospheres $\partial A$ and $\partial C$. This concept of horospherical width was first introduced in $\hyp^2$ by Fillmore \cite{Fil70}, but it has a natural extension to $\hyp^n$.

Let $K\subset \hyp^n$ be a convex body, and $A\subset \hyp^n$ a horoball. Then $\partial A$ is called a supporting horosphere to $K$, if $\partial A\cap\partial K\neq\emptyset$ and either $K\subset A$ or $K\subset \hyp^n\setminus\Int A$. For an ideal point $i$ and a convex body $K\subset \hyp^n$, there are exactly two parallel supporting horospheres orthogonal to $i$. If the horoballs determined by these supporting horospheres are $A$ and $C$, then the \emph{Fillmore-width} of $K$ at $i$ is
$$
w_F\left(K,i\right)=w_{\mathrm{hor}}\left(\partial A,\partial C\right).
$$
\begin{center}
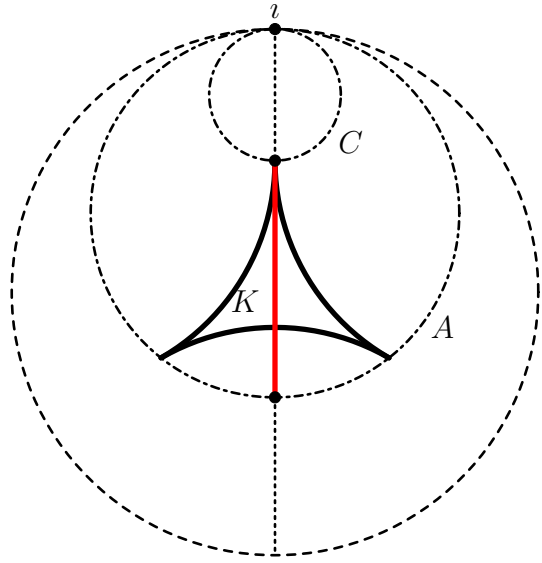
\begin{figure}
\begin{tikzpicture}[line cap=round,line join=round,>=triangle 45,x=1cm,y=1cm,scale=3.5]
\clip(-1.1,-1.1) rectangle (1.1,1.2);
\draw [line width=1pt,dashed] (0,0) circle (1cm);
\draw [shift={(-0.8660254037844384,0.5)},line width=2pt]  plot[domain=-1.0471975511965983:0,variable=\t]({1*0.8660254037844385*cos(\t r)+0*0.8660254037844385*sin(\t r)},{0*0.8660254037844385*cos(\t r)+1*0.8660254037844385*sin(\t r)});
\draw [shift={(0.866025403784439,0.5)},line width=2pt]  plot[domain=3.141592653589793:4.1887902047863905,variable=\t]({1*0.866025403784439*cos(\t r)+0*0.866025403784439*sin(\t r)},{0*0.866025403784439*cos(\t r)+1*0.866025403784439*sin(\t r)});
\draw [shift={(0,-1)},line width=2pt]  plot[domain=1.0471975511965974:2.0943951023931957,variable=\t]({1*0.8660254037844383*cos(\t r)+0*0.8660254037844383*sin(\t r)},{0*0.8660254037844383*cos(\t r)+1*0.8660254037844383*sin(\t r)});
\draw [line width=1pt,dashdotted] (0,0.3) circle (0.7cm);
\draw [line width=1pt,dashdotted] (0,0.75) circle (0.25cm);
\draw [line width=1pt,dotted] (0,1)-- (0,0.5);
\draw [line width=2pt,color=red] (0,0.5)-- (0,-0.4);
\draw [line width=1pt,dotted] (0,-0.4)-- (0,-1);
\draw (0,1) node[anchor=south] {$i$};
\draw (-0.20858126753840736,0.05) node[anchor=north west] {$K$};
\draw (0.55,-0.05063867561412254) node[anchor=north west] {$A$};
\draw (0.2,0.6529764610121744) node[anchor=north west] {$C$};
\begin{scriptsize}
\draw [fill=black] (0,-0.4) circle (.6pt);
\draw [fill=black] (0,1) circle (.6pt);
\draw [fill=black] (0,0.5) circle (.6pt);
\end{scriptsize}
\end{tikzpicture}
\caption{The Fillmore width $w_F\left(K,i\right)$ of a regular triangle $K$ at the ideal point $i$}
\end{figure}
\end{center}

\noindent{\bf Remark. }In contrast with the Santal\'o width, for the Fillmore-width it is trivially true, that for an ideal point $i$ and two convex bodies $K\subseteq L$, $w_F\left(K,i\right)\leq w_F\left(L,i\right)$. However, if $w_F\left(K,i\right)=w_F\left(L,i\right)$ for every ideal point $i$ and convex bodies $K\subseteq L$, it is not necessarily true, that $K=L$, only that the h-convex hull of $K$ and $L$ are the same. Note that in general it is not true, that for every boundary point there is a supporting horospherical strip (e.g. polytopes are counterexamples), only for h-convex bodies.

Let 
$$w_F\left(K\right)=\min\left\{w_F\left(K,i\right)\colon i\text{ is an ideal point}\right\}.$$
The following is easy to see (for example let $C_m$ be an isosceles triangle with $s$ as its base and base altitude of length $\frac{1}{m}$, or $C_m=s^{\left(\frac{1}{m}\right)}$).

\begin{prop}
Let $s=\left[x,y\right]$ be a hyperbolic line segment, and $\left\{C_m\right\}_{m\in\N}$ a sequence of convex bodies tending to $s$ with respect to the Hausdorff metric. Then,
$$
\lim_{m\to\infty}w_F\left(C_m\right)>0.
$$
\end{prop}

\subsection{Leichtweiss' width function}

For a convex body $K\subset \hyp^n$ and a hyperplane $H$ such that $H\cap K\neq\emptyset$, we say that the closed convex set $\mathcal{A}$ is the \emph{supporting hyperspherical domain} with respect to $H$, if $K\subset\mathcal{A}$, and the connected components of $\partial\mathcal{A}$ are $H_1$ and $H'_2$, where $H_1$ and $H'_2$ are components of the hyperspheres $\partial H^{\left(\varepsilon\right)}$ and $\partial H^{\left(\varepsilon'\right)}$, such that $H_1\cap K\neq\emptyset$ and $H'_2\cap K\neq\emptyset$. Then the width of the supporting hyperspherical domain $\mathcal{A}$ with respect to the hyperplane $H$ is $w_{hyp}\left(K,H\right)=\varepsilon+\varepsilon'$. The following notion of hyperspherical width was first used by Leichtweiss\cite{Lei05} in $\hyp^2$.

For a convex body $K\subset \hyp^n$ and a fixed interior point $p\in\Int K$, we define the \emph{Leichtweiss width} of $K$ with respect to a hyperplane $H$ containing $p$ as $$w_L\left(K,\left(p,H\right)\right)=w_{hyp}\left(K,H\right)$$ and we put 
$$w_L\left(K,p\right)=\min\left\{w_L\left(K,\left(p,H\right)\right)\colon H\text{ is a hyperplane containing }p\right\}.$$

\begin{center}
\begin{figure}
\begin{tikzpicture}[line cap=round,line join=round,>=triangle 45,x=1cm,y=1cm,scale=3.5]
\clip(-1.1,-1.1) rectangle (1.1,1.1);
\draw [line width=1pt,dashed] (0,0) circle (1cm);
\draw [shift={(-0.8660254037844384,0.5)},line width=2pt]  plot[domain=-1.0471975511965983:0,variable=\t]({1*0.8660254037844385*cos(\t r)+0*0.8660254037844385*sin(\t r)},{0*0.8660254037844385*cos(\t r)+1*0.8660254037844385*sin(\t r)});
\draw [shift={(0.866025403784439,0.5)},line width=2pt]  plot[domain=3.141592653589793:4.1887902047863905,variable=\t]({1*0.866025403784439*cos(\t r)+0*0.866025403784439*sin(\t r)},{0*0.866025403784439*cos(\t r)+1*0.866025403784439*sin(\t r)});
\draw [shift={(0,-1)},line width=2pt]  plot[domain=1.0471975511965974:2.0943951023931957,variable=\t]({1*0.8660254037844383*cos(\t r)+0*0.8660254037844383*sin(\t r)},{0*0.8660254037844383*cos(\t r)+1*0.8660254037844383*sin(\t r)});
\draw (-0.26,0.2) node[anchor=north west] {$K$};
\draw [line width=1pt,dotted] (-1,0)-- (1,0);
\draw [shift={(0,1.5)},line width=1pt,dashdotted]  plot[domain=4.124386376837123:5.300391583932257,variable=\t]({1*1.8027756377319966*cos(\t r)+0*1.8027756377319966*sin(\t r)},{0*1.8027756377319966*cos(\t r)+1*1.8027756377319966*sin(\t r)});
\draw [shift={(0,-0.75)},line width=1pt,dashdotted]  plot[domain=0.6435011087932844:2.498091544796509,variable=\t]({1*1.25*cos(\t r)+0*1.25*sin(\t r)},{0*1.25*cos(\t r)+1*1.25*sin(\t r)});
\draw [line width=2pt,color=red] (0,0.5)-- (0,-0.30277563773199484);
\draw (0,0) node[anchor=north west] {$p$};
\draw (0.44,0.14) node[anchor=north west] {$H$};
\begin{scriptsize}
\draw [fill=black] (0,0) circle (.6pt);
\draw [fill=black] (0,0.5) circle (.6pt);
\draw [fill=black] (0,-0.30277563773199484) circle (.6pt);
\end{scriptsize}
\end{tikzpicture}
\caption{The Leichtweiss width $w_L\left(K,\left(p,H\right)\right)$ of a regular triangle $K$ at a line $H$ that contains the inner point $p$}
\end{figure}
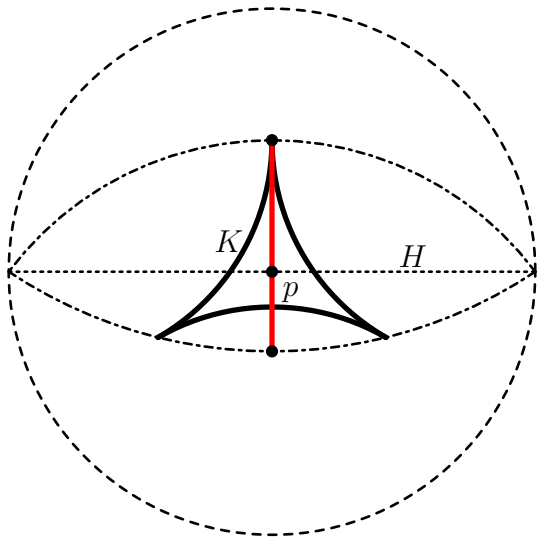
\end{center}

\noindent{\bf Remark. }We can observe, that once again, there is no supporting hyperspherical domain through all boundary points necessarily (for example polytopes). It is trivial, that for convex bodies $K_1\subseteq K_2$ and $p\in \Int K_1$, $w_L\left(K_1,\left(p,H\right)\right)\leq w_L\left(K_2,\left(p,H\right)\right)$ for any hyperplane $H$ through $p$, so $w_L\left(K_1,p\right)\leq w_L\left(K_2,p\right)$. Note that $w_L\left(K,p\right)$ depends on the choice of $p$, for instance, if $K$ is an equilateral triangle in $\hyp^2$ and $\ell$ is one of its altitudes, then $w_L\left(K,p\right)$ gets greater as $p\in\ell$ gets closer to the vertex, and gets smaller (approximating the altitude of the triangle) as $p$ gets closer to the midpoint of the opposite side.

\subsection{Planar alternatives}

There are two more recent approaches to define hyperbolic width on the plane $\hyp^2$. The first one is by Jer\'onimo-Castro and Jimenez-Lopez\cite{JCJL17}. 
If $z\in\partial K$ is a boundary point and $\ell_z$ is a supporting line through $z$, we call a $z'\in\partial K$ an \emph{opposite point} to $z$ with respect to $\ell_z$, if there is a supporting line $\ell_{z'}$ such that $\ell_z$ and $\ell_{z'}$ have the same alternating angles with the geodesic segment $\left[z,z'\right]$. By definition, the \emph{Jer\'onimo-Castro--Jimenez-Lopez} width of $K$ at $\left(z,\ell_z,z'\right)$ is
$$
w_{JCJL}\left(K,\left(z,\ell_z,z'\right)\right)=d\left(\ell_z,\ell_{z'}\right).
$$
 Note that for a convex body $K$ and a supporting line $\ell_z$ through $z\in\partial K$ the opposite point $z'$ is not unique in general, but if $K$ is h-convex, then it is unique (see \cite{JCJL17}).

Another more recent planar width function is defined by G. Horv\'ath\cite{Hor21}. For a convex body $K\subset \hyp^2$ and an ideal point $i$, we define the \emph{G. Horv\'ath width} as
$$
w_{GH}\left(K,i\right)=\sup\left\{w_{\mathrm{hyp}}\left(\ell\right)\colon\ell\cap\Int K\neq\emptyset\text{ and }i\text{ is an ideal point of }\ell\right\}.
$$

We put
$$w_{JCJL}(K)=\min\{w_{JCJL}(K,(z,\ell_z,z'))\colon z,z'\in\partial K\text{ are opposite with respect to }\ell_z\}$$
and
$$w_{GH}(K)=\min\{w_{GH}(K,i)\colon i\text{ is an ideal point}\}.$$

\subsection{A new concept of width}

We will extend the width function $w_{hyp}$ introduced by Leichtweiss, and hence we will refer to this new concept of width as the extended Leichtweiss width. Let $K\subset \hyp^n$ be a convex body and let $H$ be a hyperplane intersecting $K$. We extend the notion of supporting hyperspherical domain $\mathcal{A}$ the following way. If $H$ is a supporting hyperplane, then there is a positive $\varepsilon$, such that $K\subset H^{\varepsilon}$, and $\partial K$ intersects one of the connected components of $\partial H^{\varepsilon}$ (say $H_1$). Then, the supporting hyperspherical domain $\mathcal{A}$ will be bounded by $H$ and $H_1$, and the corresponding width is $w_{hyp}\left(H\right)=\varepsilon$ in this case. Now let $H$ be an arbitrary hyperplane intersecting $K$. Then the \emph{width} of $K$ with respect to the hyperplane $H$ is $w\left(K,H\right)=w_{hyp}\left(K,H\right)$.

We define the minimal width of $K$ as
$$w(K)=\min\{w(K,H)\colon H\text{ is a hyperplane intersecting }K\}$$

\noindent\textbf{Remark. }It is clear, that for $n=2$ we have
$$w_{GH}(K,i)=\max\{w(K,\ell)\colon i\text{ is an ideal point of }\ell\}$$


\section{Continuity of width functions}
\label{cont}

It is natural to ask whether the above defined notions of width are continuous. To answer this we first have to define what we mean by continuity for these objects. In the following we will introduce the necessary definitions and then the rest of this section concerns the proof of Theorem \ref{contthm}, i.e. that the width functions are indeed continuous in every case. 
We will work in the Poincar\'e disk model.

We say that {\it a sequence of ideal points $i_m$ converges to an ideal point} $i$ if it does so in the usual topology of the sphere $\partial B^n$. This convergence is well-defined since the identification of $\hyp^n$ with $\mathrm{int}B^n$ is a homeomorphism and so if we choose another Poincar\'e disk model for $\hyp^n$, the convergent ideal point sequences remain the same.

Let $H_m,H$ be hyperplanes with ideal point sets $I_m,I\subset\partial B^n$ respectively and let $c_m,c$ and $r_m,r$ be the centres and radii of $I_m,I$ in the sphere $\partial B^n$ respectively. Note that the centre-radius pair $(-c,\pi-r)$ yields the same hyperplane as the pair $(c,r)$ (where we think of $B^n$ as the unit ball in $\R^n$ and $-c$ denotes the reflection of $c$ to the origin). Then {\it the sequence $H_m$ converges to $H$} if $\{(c_m,r_m),(-c_m,\pi-r_m)\}\to\{(c,r),(-c,\pi-r)\}$.
Observe that although the centres and radii of $I_m$ and $I$ in $\partial B^n$ can change when changing the model for $\hyp^n$, the notion of convergence is again independent of the model we choose.

If $\ell_m$ and $\ell$ are lines with ideal points $i_m,j_m$ and $i,j$ respectively, then {\it the sequence $\ell_m$ converges to $\ell$} if $\{i_m,j_m\}\to\{i,j\}$.

\noindent\textbf{Remark. }With the above definitions the space of ideal points is $\partial B^n$ which is compact. The space of hyperplanes is
$$\mathcal{H}:=\bigslant{(\partial B^n\times(0,\pi))}{(c,r)\sim(-c,\pi-r)}$$
and the space of lines is
$$\mathcal{L}:=\bigslant{\{(i,j)\in\partial B^n\times\partial B^n\colon i\ne j\}}{(i,j)\sim(j,i)}$$
and both are locally compact.

We will use the following basic topological lemma.

\begin{lemma}
\label{toplemma}
Let $X$ be a locally compact and $Y$ a compact topological space, let $f\colon X\times Y\to\R$ be a continuous function and for each $x\in X$ fix a non-empty closed subset $Z_x\subset Y$. Then the function
$$F\colon X\to\R;~x\mapsto\max_{y\in Z_x}f(x,y)$$
is also continuous.
\end{lemma}

We will usually use this lemma with $Z_x=Y$ for all $x$.

\proof
Let $x_m\to x$ be a convergent sequence in $X$ and fix an $\varepsilon>0$; we want to prove that if $m$ is sufficiently large, then
$$|\max_{y\in Z_x}f(x,y)-\max_{y\in Z_{x_m}}f(x_m,y)|<\varepsilon.$$
Since $X$ is locally compact, we may assume that $x$ and $x_m$ are all in a compact subspace $C\subset X$. Now $C\times Y$ is compact, hence there are open subsets $U\subset C$ and $V_1,\ldots,V_n\subset Y$ such that $x\in U$ and for all $i$ if $a,b\in U\times V_i$, then $|f(a)-f(b)|<\varepsilon$. Let $y_m,y\in Y$ be points such that
$$f(x_m,y_m)=\max_{Z_{x_m}}f(x_m,\cdot)\quad\text{and}\quad f(x,y)=\max_{Z_x}f(x,\cdot).$$
If $m$ is sufficiently large, then $x_m\in U$ which implies $|f(x,y_m)-f(x_m,y_m)|<\varepsilon$ and $|f(x_m,y)-f(x,y)|<\varepsilon$. Now if $f(x,y)\ge f(x_m,y_m)$, then we get
$$f(x_m,y_m)-\varepsilon\le f(x,y)-\varepsilon<f(x_m,y)\le f(x_m,y_m)$$
thus $|f(x,y)-f(x_m,y_m)|<\varepsilon$; similarly if $f(x_m,y_m)\ge f(x,y)$, we get
$$f(x,y)-\varepsilon\le f(x_m,y_m)-\varepsilon<f(x,y_m)\le f(x,y)$$
so again $|f(x,y)-f(x_m,y_m)|<\varepsilon$. This is what we wanted to prove.
\endproof

\begin{lemma}
\label{santalolemma}
Let $K\subset \hyp^n$ be a convex body, $\ell_m\to\ell$ a convergent sequence of lines with ideal points $i_m\to i$, $j_m\to j$ and let $H_m,H'_m$ (resp. $H,H'$) be the supporting hyperplanes of $K$ orthogonal to $\ell_m$ (resp. $\ell$) such that $H_m\cap\ell_m\in[i_m,H'_m\cap\ell_m]$ and $H\cap\ell\in[i,H'\cap\ell]$. Then with the notations $\{a_m\}=H_m\cap\ell_m$, $\{a'_m\}=H'_m\cap\ell_m$ and $\{a\}=H\cap\ell$, $\{a'\}=H'\cap\ell$ we have $d_{\hyp^n}(a_m,a'_m)\to d_{\hyp^n}(a,a')$.
\end{lemma}

\proof
Let $I\subset\partial B^n$ be the set of ideal points of $H$ and let $\mathcal{L}_H$ be the space of those lines in $\hyp^n$ whose two ideal points are in different components of $\partial B^n\setminus I$. Note that $\mathcal{L}_H$ is locally compact, $\ell\in \mathcal{L}_H$ and we can assume that $\ell_m\in \mathcal{L}_H$ for all $m$. Define
$$f_H\colon \mathcal{L}_H\times K\to\R$$
in a point $(k,x)$ in the following way: let $y$ be the intersection point of the line $k$ with $H$ and let $z$ be the point on $k$ closest to $x$ (that is, the intersection of $k$ with the hyperplane through $x$ orthogonal to $k$); then we put $f_H(k,x):=d(y,z)$ if $z$ is in the component of $\hyp^n\setminus H$ whose closure contains the ideal point $i$ and $f_H(k,x):=-d(y,z)$ otherwise. Now clearly
$$d_{\hyp^n}(a,a')=\max_{x\in K}f_H(\ell,x)-\min_{x\in K}f_H(\ell,x)\quad\text{and}\quad d_{\hyp^n}(a_m,a'_m)=\max_{x\in K}f_H(\ell_m,x)-\min_{x\in K}f_H(\ell_m,x),$$
hence we only need to prove that $\max f_H(\ell_m,\cdot)\to\max f_H(\ell,\cdot)$ and $\min f_H(\ell_m,\cdot)\to\min f_H(\ell,\cdot)$. This will immediately follow from Lemma \ref{toplemma} as soon as we prove that $f_H$ is continuous.

Let $(k_m,x_m)\to(k,x)$ be a convergent sequence of lines in $\mathcal{L}_H$ and points in $K$, moreover, let $y,z$ be the points on $k$ whose distance defines the value of $f_H(k,x)$ as above and let $y_m,z_m$ be their analogues on $k_m$ defining $f_H(k_m,x_m)$. It suffices to prove the convergences $y_m\to y$ and $z_m\to z$. For the proof of $y_m\to y$ we use the Beltrami--Cayley--Klein model where hyperbolic lines and hyperplanes are represented as the intersections of $\mathrm{int}B^n$ with Euclidean lines and hyperplanes respectively. Clearly the convergence of ideal points (and so of lines) is the same as in the Poincar\'e disk model and now the points $y_m$ and $y$ are realised as intersections of a Euclidean hyperplane with Euclidean line segments whose endpoints converge, hence the convergence $y_m\to y$ follows. For the proof of $z_m\to z$ let $z'_m$ denote the point on $k_m$ closest to $x$. Choosing a Beltrami--Cayley--Klein model where $x$ is the centre of $B^n$ yields that $z'_m$ is represented as the intersection of $k_m$ with the Euclidean hyperplane through $x$ and orthogonal to $k_m$, hence we get $z'_m\to z$. Since $x_m\to x$, for any $\varepsilon>0$ we have $|d(x,z'_m)-d(x_m,z'_m)|<\varepsilon$ and $|d(k_m,x)-d(k_m,x_m)|<\varepsilon$ for sufficiently large indices $m$, which implies
$$|d(x_m,z_m)-d(x_m,z'_m)|\le|d(x_m,z_m)-d(x,z'_m)|+|d(x,z'_m)-d(x_m,z'_m)|<2\varepsilon.$$
Thus $|d(z_m,z'_m)|\to0$ and so $z_m\to z$ which finishes the proof.
\endproof

\begin{prop}
The Santal\'o width $w_S(K,(z,H_z))$ is a continuous function of $(K,(z,H_z))$.
\end{prop}

\proof
Let $K_m\to K$ be a convergent sequence of convex bodies and $H_m$ supporting hyperplanes through $z_m\in\partial K_m$ such that $(z_m,H_m)\to(z,H_z)$. If $\ell_m$ and $\ell$ denote the lines through $z_m$ and $z$ orthogonal to $H_m$ and $H$ respectively, then it is not hard to see that $\ell_m\to\ell$. Now the previous lemma implies that we have $w_S(K_m,(z_m,H_m))\to w_S(K,(z,H_z))$, hence $w_S$ is indeed continuous.
\endproof 

\begin{lemma}
\label{fillmorelemma}
If $K\subset \hyp^n$ is a convex body, $i$ an ideal point and $i_m$ a sequence of ideal points with $i_m\to i$, then we have $w_F(K,i_m)\to w_F(K,i)$.
\end{lemma}

\proof
Let $o$ be the centre of the disk $B^n$ (where again $\mathrm{int}B^n$ is identified with $\hyp^n$). We define the value of
$$h_o\colon\partial B^n\times K\to\R$$
at a point $(j,x)$ as follows: let $y$ be the intersection of the horosphere at $j$ through $x$ with the geodesic line through $o$ tending to $j$ oriented positive from $o$ in the direction of $j$ and negative in the direction of $-j$; then $h_o(j,x)$ is the oriented distance of $y$ from $o$. Now we have
$$w_F(K,j)=\max_{x\in K}h_o(j,x)-\min_{x\in K}h_o(j,x),$$
hence we only need to prove that $\max h_o(i_m,\cdot)$ and $\min h_o(i_m,\cdot)$ converge to $\max h_o(i,\cdot)$ and $\min h_o(i,\cdot)$ respectively. Since $h_o$ is clearly continuous, we can now apply Lemma \ref{toplemma} with $X=\partial B^n$, $Y=K$ and $f=h_o$ for the maximum and $f=-h_o$ for the minimum.
\endproof

\begin{prop}
\label{fillmorecont}
The Fillmore width $w_F(K,i)$ is a continuous function of $(K,i)$.
\end{prop}

\proof
Let $(K_m,i_m)\to(K,i)$ be a convergent sequence of pairs of convex bodies and ideal points and fix an $\varepsilon>0$. If $m$ is large enough, then $K_m$ and $K$ are contained in the convex bodies $K^{(\varepsilon)}$ and $K_m^{(\varepsilon)}$ respectively, hence for any ideal point $j$ we have
$$w_F(K,j)-2\varepsilon\le w_F(K_m^{(\varepsilon)},j)-2\varepsilon\le w_F(K_m,j)\le w_F(K^{(\varepsilon)},j)\le w_F(K,j)+2\varepsilon$$
Now the previous lemma implies that for $m$ sufficiently large we have $|w_F(K,i_m)-w_F(K,i)|<\varepsilon$, thus we obtain
$$|w_F(K_m,i_m)-w_F(K,i)|\le|w_F(K_m,i_m)-w_F(K,i_m)|+|w_F(K,i_m)-w_F(K,i)|<5\varepsilon$$
for all suitably large indices $m$, which finishes the proof.
\endproof

\begin{lemma}
\label{leichtweisslemma}
If $K\subset \hyp^n$ is a convex body, $H$ a hyperplane intersecting $K$ and $H_m$ a sequence of hyperplanes intersecting $K$ such that $H_m\to H$, then we have $w(K,H_m)\to w(K,H)$.
\end{lemma}

\proof
Let $I\subset\partial B^n$ be the set of ideal points of $H$ in a Poincar\'e disk model and let $c,-c$ be the two (antipodal) centre points of $I$ on the sphere. Denote by $\mathcal{H}_H$ the space of those hyperplanes $H'\subset \hyp^n$ whose analogous two centre points are in different components of $\partial B^n\setminus I$. Observe that $\mathcal{H}_H$ is locally compact, $H\in \mathcal{H}_H$ and we can also assume that $\mathcal{H}_H$ contains $H_m$ for all $m$. We define the value of
$$g_H\colon \mathcal{H}_H\times K\to\R$$
in a hyperplane $H'$ and point $x$ as the distance of $x$ and $H'$ with positive sign if $x$ is in the half-space whose closure contains $-c$ and with negative sign otherwise. Then we have
$$w(K,H')=\max_{x\in K}g_H(H',x)-\min_{x\in K}g_H(H',x)$$
for any hyperplane $H'$ intersecting $K$, thus it suffices to prove the convergences $\max g_H(H_m,\cdot)\to\max g_H(H,\cdot)$ and $\min g_H(H_m,\cdot)\to\min g_H(H,\cdot)$. Now we can apply Lemma \ref{toplemma} with $f=g_H$ for the maximum and $f=-g_H$ for the minimum provided that $g_H$ is continuous.

If we have a convergent sequence of hyperplanes $H'_m\to H'$ in $\mathcal{H}_H$ and a convergent sequence of points $x_m\to x$ in $K$, then we need that the distances $d(H'_m,x_m)$ converge to $d(H',x)$. This follows analogously to the way we obtained the convergence $z_m\to z$ in the proof of Lemma \ref{santalolemma} so we omit its proof here.
\endproof

Now a word-by-word analogue of the proof of Proposition \ref{fillmorecont} yields the following.

\begin{prop}
\label{wKHcont}
The extended Leichtweiss width $w(K,H)$ is a continuous function of $(K,H)$.
\end{prop}

As a special case of this we also get the following.

\begin{prop}
The Leichtweiss width $w_L(K,(p,H_p))$ is a continuous function of $(K,(p,H_p))$.
\end{prop}

\begin{lemma}
Let $z_m\to z$ and $z'_m\to z'$ be convergent sequences of points in $\hyp^2$ and let $\ell_m\to\ell_z$ and $\ell'_m\to\ell_{z'}$ be convergent sequences of lines such that $z_m\in\ell_m$, $z\in\ell_z$, $z'_m\in\ell'_m$ and $z'\in\ell_{z'}$, the sequence $\ell_m$ converges to $\ell_z$ and  $\ell_m$ and $\ell'_m$ (resp. $\ell_z$ and $\ell_{z'}$) have equal alternating angles with the segment $[z_m,z'_m]$ (resp. $[z,z']$). Then we have $d(\ell_m,\ell'_m)\to d(\ell_z,\ell_{z'})$.
\end{lemma}

\proof
Since the distance of lines is continuous, it is sufficient to prove that $\ell'_m\to\ell_{z'}$. Let $\varphi_m$ denote the angle of the line $\ell_m$ and the segment $[z_m,z'_m]$ which is oriented positive from the line to the segment; similarly let $\varphi$ be the positive angle from $\ell_z$ to $[z,z']$. Observe that $\varphi_m$ tends to $\varphi$ and the positive angle from the line $\ell'_m$ (resp. $\ell_{z'}$) to the segment $[z'_m,z_m]$ (resp. $[z',z]$) is again $\varphi_m$ (resp. $\varphi$). Now if $\ell''_m$ is the line through $z'_m$ which has positive angle $\varphi$ with $[z'_m,z_m]$, then it is not hard to see the convergence $\ell''_m\to\ell'_z$. Let $i'_m$, $j'_m$ and $i''_m$, $j''_m$ be the ideal points of $\ell'_m$ and $\ell''_m$ respectively such that $\varphi_m$ is the angle $\angle(i'_m,z'_m,z_m)$ and $\varphi=\angle(i''_m,z'_m,z_m)$. Now the angle of $\ell'_m$ and $\ell''_m$ tends to $0$, hence if $B^2$ denotes the disk of a Poincar\'e disk model we use, the sequence of the (Euclidean) distances $d_{\partial B^2}(i'_m,i''_m)$ and $d_{\partial B^2}(j''_m,j'_m)$ of the corresponding ideal points both tend to zero. This implies that $\ell'_m$ also converges to $\ell_{z'}$ which concludes the proof.
\endproof

\begin{prop}
The Jer\'onimo-Castro--Jimenez-Lopez width $w_{JCJL}(K,(z,\ell_z,z'))$ is a continuous function of $(K,(z,\ell_z,z'))$. 
\end{prop}

\proof
Let $K_m\to K$ be a convergent sequence of convex bodies and let $z_m,z'_m\in\partial K_m$ two sequences of points converging to $z,z'\in\partial K$ respectively. Let $\ell_m\to\ell_z$ be a convergent sequence of lines such that $\ell_m$ (resp. $\ell_z$) is a supporting line of $K_m$ at $z_m$ (resp. of $K$ at $z$) and suppose that the point $z'_m$ (resp. $z'$) is an opposite point to $z_m$ (resp. $z$) with respect to $\ell_m$ (resp. $\ell$). Denote by $\ell'_m$ the supporting line through $z'_m$ such that $\ell'_m$ and $\ell_m$ have equal alternating angles with $[z_m,z'_m]$ and let $\ell_{z'}$ be the analogous supporting line through $z'$. Now since we have $w_{JCJL}(K_m,(z_m,\ell_m,z'_m))=d(\ell_m,\ell'_m)$ and $w_{JCJL}(K,(z,\ell_z,z'))=d(\ell_z,\ell_{z'})$, the previous lemma yields that $w_{JCJL}(K_m,(z_m,\ell_m,z'_m))\to w_{JCJL}(K,(z,\ell_z,z'))$. This shows that $w_{JCJL}(K,(z,\ell_z,z'))$ is indeed continuous.
\endproof

\begin{lemma}
If $K\subset \hyp^2$ is a convex body and $i_m\to i$ is a convergent sequence of ideal points, then we have $w_{GH}(K,i_m)\to w_{GH}(K,i)$.
\end{lemma}

\proof
Recall that we have $w_{GH}(K,j)=\max\{w(K,\ell)\colon j\text{ is an ideal point of }\ell\}$. Now putting $v(j,x):=w(K,\ell_{j,x})$ where $\ell_{j,x}$ is the line with ideal point $j$ through $x$, Proposition \ref{wKHcont} implies that $v$ is continuous and by definition we get
$$w_{GH}(K,j)=\max_{x\in K}v(j,x).$$
Now we can apply Lemma \ref{toplemma} with $X=\partial B^n$ and $Y=K$ which yields that $w_{GH}(K,j)$ is continuous in the variable $j$, thus our claim follows.
\endproof

Now again a direct analogue of the proof of Proposition \ref{fillmorecont} implies the following.

\begin{prop}
The G. Horv\'ath width $w_{GH}(K,i)$ is a continuous function of $(K,i)$.
\end{prop}

\noindent\textbf{Remark. }It is worthwhile to note that in the proofs above we never used the convexity of the bodies $K$ and $K_m$, hence the support functions $f_H$, $g_H$ and $h_o$ defined in the proofs of the lemmas \ref{santalolemma}, \ref{fillmorelemma} and \ref{leichtweisslemma} can be used to extend these notions of width to any compact subset of $\hyp^n$ (and with these definitions the width of any compact subset will be the width of its convex hull). The width functions will then be continuous in this extended sense too.

We finish this section with showing that the minimal width functions defined by all of the above widths are continuous as well.

\begin{theorem}
\label{minwcont}
If $K_m\to K$ is a convergent sequence of convex bodies in $\hyp^n$, then
\begin{description}
\item[(i)] $w_S(K_m)\to w_S(K)$;
\item[(ii)] $w_F(K_m)\to w_F(K)$;
\item[(iii)] if $p_m\in K_m$ and $p\in K$ are interior points and $p_m\to p$, then $w_L(K_m,p_m)\to w_L(K,p)$;
\item[(iv)] if $n=2$, then $w_{JCJL}(K_m)\to w_{JCJL}(K)$;
\item[(v)] if $n=2$, then $w_{GH}(K_m)\to w_{GH}(K)$;
\item[(vi)] $w(K_m)\to w(K)$.
\end{description}
\end{theorem}

\proof
As we saw above, all of these widths are continuous as functions of $K$ and a second object, which is either an ideal point or a hyperplane (or line) intersecting or supporting $K$ occasionally with fixed additional points. We want to apply Lemma \ref{toplemma} to the negative of these width functions with $X$ being the space of compact subsets and $Y$ the space of the objects in the second variable in all cases. The Blaschke Selection Theorem \ref{BlaschkeSel} implies that this $X$ is locally compact, hence if $Y$ is the space of ideal points (which is compact), i.e. for $w_F$ and $w_{GH}$, then the lemma may indeed be applicable. To be able to apply the lemma in the other cases, note that it is sufficient to use it with $X$ being only the space of compact subsets of $K^{(\varepsilon)}$ for some $\varepsilon>0$. Then we can consider only those hyperplanes which intersect $K^{(\varepsilon)}$, which form a closed subspace
$$\mathcal{H}_{K,\varepsilon}\subset\bigslant{(\partial B^n\times[\delta,\pi-\delta])}{(c,r)\sim(-c,\pi-r)}\subset\mathcal{H}$$
for some $\delta>0$, hence it is compact. Now for the width functions $w_L$ and $w$ we can put $Y:=\mathcal{H}_{K,\varepsilon}$ and for $w_S$ and $w_{JCJL}$ the product spaces $\mathcal{H}_{K,\varepsilon}\times K^{(\varepsilon)}$ and $K^{(\varepsilon)}\times\mathcal{H}_{K,\varepsilon}\times K^{(\varepsilon)}$ respectively. Then we can define for each convex body $K'$ the subspace $Z_{K'}\subset Y$ as the set of those objects where the appropriate width function is defined, which is closed in all cases. Hence Lemma \ref{toplemma} can be applied and our statement follows.
\endproof

This finishes the proof of Theorem \ref{contthm}.

\section{Width and diameter}
\label{diam}

In what follows we will consider the maximal widths defined by the width functions $w_S$, $w_F$, $w_L$, $w_{JCJL}$, $w_{GH}$ and $w$ and we shall see that this always gives the diameter.

\begin{lemma}\label{diam_Santalo}
For a convex body $K\subset \hyp^n$,
$$
\diam K=\max\left\{w_S\left(K,\left(z,H_z\right)\right)\colon z\in\partial K\text{ and }H_z\text{ is a supporting hyperplane of }K\text{ at }z\right\}.
$$
\end{lemma}
\proof
Let $z\in\partial K$ be an arbitrary boundary point and $H_z$ any supporting hyperplane to $K$, such that $z\in H_z$. If $\ell$ is the unique hyperbolic line orthogonal to $H_z$ through the point $z$, and $H'_z$ is the other supporting hyperplane to $K$ orthogonal to $\ell$, intersecting $\ell$ at the point $z'$, then there is some point $z_0\in H'_z\cap\partial K$. Because of the orthogonality, and from the definition of the diameter,
$$
w_S\left(K,\left(z,H_z\right)\right)=d\left(z,z'\right)\leq d\left(z,z_0\right)\leq\diam K
$$
follows. On the other hand, if $\left[z,z'\right]$ is a diametral chord in $K$, then the hyperplanes $H_z$ and $H_{z'}$ through the points $z$ and $z'$ respectively and both orthogonal to the geodesic line connecting $z$ and $z'$ are supporting $K$ at $z$ and $z'$, respectively. Otherwise, if there is a point, say $z_0\in K$ in the open half-space bounded by $H_{z'}$ not containing $H_z$, then
$$
d\left(z,z_0\right)>d\left(z,z'\right)=\diam K
$$
is a contradiction. Therefore, in this case
$$
w_S\left(K,\left(z,H_z\right)\right)=w_S\left(K,\left(z',H_{z'}\right)\right)=d\left(z,z'\right)=\diam K.
$$
This concludes the proof.
\endproof

\begin{lemma}
Let $K\subset \hyp^n$. Then,
$$
\diam K=\max\left\{w_F\left(K,i\right)\colon i\text{ is an ideal point}\right\}.
$$
\end{lemma}
\proof
Let $i$ be any ideal point, and let $\partial A$ and $\partial C$ be the two parallel supporting horospheres to $K$ at the ideal point $i$. Then, there is a point $x\in\partial A\cap \partial K$ and $y\in\partial C\cap\partial K$. Let $z$ be the intersection of the horosphere $\partial C$ and the line through the ideal point $i$ and $x$, orthogonal to both of these supporting horospheres. Then,
$$
\diam K\geq d\left(x,y\right)\geq d\left(x,z\right)=w_F\left(K,i\right).
$$
On the other hand, if $\left[x,y\right]\subset K$ is a diametral chord, and $i$ is one of the ideal points determined by the geodesic line $\ell$ through $x$ and $y$, then if $\partial A$ is the orthogonal horosphere to $i$ through $x$ and $\partial C$ is the orthogonal horosphere to $i$ containing $y$, $\partial A$ and $\partial C$ need to be supporting horospheres of $K$. Otherwise, if there is a point say $z\in\Int A\cap K$, where we assume $A\subset C$, then
$$
d\left(y,z\right)\geq d\left(x_0,z\right)+d\left(x_0,y_0\right)>d\left(x_0,y_0\right)=\diam K
$$
is a contradiction, where $x_0$ is the intersection point of the line $\ell_z$ through $i$ and $z$ with the horosphere $\partial A$ and $y_0$ is the intersection of $\ell_z$ and $\partial C$.
\endproof

\begin{lemma}\label{diam_Leichtweiss}
Let $K\subset \hyp^n$ be a convex body and $p\in\Int K$ a fixed interior point. Then,
$$
\diam K=\max\left\{w_L\left(K,\left(p,H\right)\right)\colon H\text{ is a hyperplane containing }p\right\}.
$$
\end{lemma}
\proof
Let $H$ be any hyperplane containing $p$, and $\mathcal{A}$ be the hyperspherical domain with respect to $H$ with boundary components $H_1$ and $H_2$. Then, there are some points $x_1\in H_1\cap\partial K$ and $x_2\in H_2\cap\partial K$. Let $x$ be the intersection point of the geodesic segment $\left[x_1,x_2\right]$ and the hyperplane $H$, and let $y_1$ and $y_2$ be the orthogonal projection of $x_1$ and $x_2$ onto $H$, respectively. Then,
$$
\diam K\geq d\left(x_1,x_2\right)=d\left(x_1,x\right)+d\left(x_2,x\right)\geq d\left(x_1,y_1\right)+d\left(x_2,y_2\right)=w_L\left(K,\left(p,H\right)\right).
$$
On the other hand, if $\left[x_1,x_2\right]\subset K$ is a diametral chord, then the hyperplane $H$ containing $p$ and orthogonal to $\left[x_1,x_2\right]$ determines a hyperspherical domain $\mathcal{A}$ of boundary components $H_1\ni x_1$ and $H_2\ni x_2$. This has to be a supporting hyperspherical domain, otherwise if there is a point $y\in\Int K\setminus\mathcal{A}$ (say in $H_1^-$) and $\ell$ denotes the line through $y$ and orthogonal to $H$, while $y_1$ and $y_2$ are the intersection points $\ell\cap H_1$ and $\ell\cap H_2$, then
$$
d\left(y,x_2\right)\geq d\left(y,y_2\right)=d\left(y,y_1\right)+d\left(y_1,y_2\right)=d\left(y,y_1\right)+d\left(x_1,x_2\right)>d\left(x_1,x_2\right)=\diam K
$$
is a contradiction.
\endproof

\begin{lemma}
For a convex body $K\subset \hyp^n$,
$$
\diam K=\max\left\{w\left(K,H\right)\colon H\text{ is a hyperplane intersecting }K\right\}.
$$
\end{lemma}
\proof
If there is an interior point $p\in\Int K\cap H$, then we know from Lemma~\ref{diam_Leichtweiss} that $w\left(K\right)\leq\diam K$. Similarly, if $H$ is a supporting hyperplane to $K$ and $p\in\partial K\cap H$, $\ell$ is the geodesic line through $p$ orthogonal to $H$ then for a sequence of interior points $\left\{p_m\right\}_{m\in\N}\subset\Int K\cap\ell$ such that $\lim_{m\to\infty}p_m=p$ and hyperplanes $H_m$ orthogonal to $\ell$ through the point $p_m$, from Proposition~\ref{wKHcont}, we can deduce
$$
w\left(K,H\right)\leq\limsup_{m\in\N} w\left(K,H_m\right)\leq\diam K.
$$
On the other hand, if $\left[z_1,z_2\right]\subset K$ is a diametral chord, $\ell$ denotes the geodesic line through $z_1$ and $z_2$, then as we have seen in Lemma~\ref{diam_Leichtweiss}, that for any hyperplane $H$ orthogonal to $\ell$ and intersecting $K$ in some interior point, $w\left(K,H\right)=\diam K$. Moreover, if $H$ is a supporting hyperplane orthogonal to $\ell$ to $K$ (say at $z_1$), then for some interior points $\left\{p_m\right\}_{m\in\N}\subset\Int\left[z_1,z_2\right]$ such that $\lim_{m\to\infty}p_m=z_1$ and orthogonal hyperplanes $H_m$ to $\ell$ through $p_m$, as $w\left(K,H_m\right)=\diam K$ for all $m\in\N$, Proposition~\ref{wKHcont} implies
$$
w\left(K,H\right)=\lim_{m\to\infty}w\left(K,H_m\right)=\diam K.
$$
\endproof

\noindent{\bf Remark. }
For $n=2$ we have
$$
\max\left\{w_{GH}\left(K,i\right)\colon i\text{ ideal point}\right\}=\max\left\{w\left(K,\ell\right)\colon\ell\text{ is a line intersecting }K\right\},
$$
so in this case the maximum of the G. Horv\'ath width also equals to the diameter by the above lemma.

G. Horv\'ath proved\cite{Hor21} the following comparison of the Santal\'o width and the Jer\'onimo-Castro--Jimenez-Lopez width on the hyperbolic plane.
\begin{lemma}\label{wJCJLwS}
Let $K\subset \hyp^2$ be a convex body, let $z\in\partial K$ be any boundary point, $\ell_z$ a supporting hyperplane at $z$ to $K$ and $z'$ an opposite point to $z$ with respect to $\ell_z$. Then
$$
w_{JCJL}\left(K,\left(z,\ell_z,z'\right)\right)\geq w_S\left(K,\left(z,\ell_z\right)\right).
$$
\end{lemma}

This and Lemma~\ref{diam_Santalo} implies the following.

\begin{corollary}\label{diam_JCJL}
Let $K\subset \hyp^2$ a convex body. Then,
$$
\diam K=\max\left\{w_{JCJL}\left(K,\left(z,\ell_z,z'\right)\right)\colon z,z'\in\partial K\text{ are opposite with respect to }\ell_z\right\}.
$$
\end{corollary}


\noindent\textbf{Remark. }A similar reasoning as in the proof of Theorem \ref{minwcont} also yields the statement of Lemma \ref{limit}, i.e. that the diameter (which is the maximal width instead of the minimal) is continuous.

\section{Bodies of constant width in the hyperbolic space}
\label{constwidth}

We call a $K\subseteq \hyp^n$ set \emph{complete} if for each $y\in \hyp^n\setminus K$, $\diam\left(K\cup\left\{y\right\}\right)>\diam\left(K\right)$. Readily, any complete set is closed. A complete set $\widetilde{K}$ is called a \emph{completion} of $K$ if $K\subseteq\widetilde{K}$ and they have the same diameter.
The Zorn lemma already yields the existence of completions (which is in general not unique, see Dekster\cite{Dek95}), but as Theorem~\ref{scott} will show, we can also give a construction of a completion, such that the completion has the same circumscribed ball as the original body.

For the notion of bodies of constant width, we use Dekster's definition (see\cite{Dek95}). A convex body $K\subset \hyp^n$ is \emph{of constant width} $D>0$, if for any boundary point $x\in\partial K$ and any outer normal unit vector $v\in T_x$, there is a boundary point $y\in\partial K$ such that $\left[x,y\right]\subset K$ is of length $D$ and $-v$ points along the geodesic segment $\left[x,y\right]$. We give the following characterization for bodies of constant width in $\hyp^n$.

\medskip\noindent\textbf{Theorem \ref{completeness}}\textit{
Let
$K\subset \hyp^n$ be a convex body of diameter $D>0$. The following are equivalent:
\begin{description}
\item[(i)] $K$ is complete;
\item[(ii)] $K=\bigcap_{x\in K}B\left(x,D\right)$;
\item[(iii)] $K$ is of constant width $D$;
\item[(iv)] $w_S\left(K,\left(z,H_z\right)\right)=D$ for all $z\in\partial K$ and $H_z$ supporting hyperplanes to $K$ at $z$;
\item[(v)] $w_F\left(K,i\right)=D$ for all ideal points $i$;
\item[(vi)] $w_L\left(K,\left(p,H_p\right)\right)=D$ for all hyperplanes through a fixed inner point $p\in\Int K$;
\item[(vii)] $w\left(K,H\right)=D$ for all hyperplanes $H$ intersecting $K$.
\end{description}}
\proof
For the equivalence of (i), (ii), (iii) and (v), see B\"or\"oczky, Sagmeister\cite{BoS22}.

Assume that $K$ is of constant width $D$ and let $z\in\partial K$ be any boundary point and $H_z$ any supporting hyperplane at $z$ to $K$. Let $\ell$ be the geodesic line containing $z$ orthogonal to $H_z$, $z'$ will be the intersection point of $\ell$ and $\partial K\setminus\left\{z\right\}$, while $v\in T_z$ is the outer normal unit vector, such that $-v$ points along the geodesic segment $\left[z,z'\right]$. Then, since $K$ is of constant width, the geodesic segment $\left[z,z'\right]$ is a diametral chord, so as we have seen in Lemma~\ref{diam_Santalo}, the hyperplane $H_{z'}$ containing $z'$ and orthogonal to $\ell$ supports $K$, therefore $w_S\left(K,\left(z,H_z\right)\right)=D$, hence (iii) implies (iv).

Now we assume that $w_S\left(K,\left(z,H_z\right)\right)=D$ for all boundary point $z\in\partial K$ and $H_z$ supporting hyperplane of $K$ at $z$. Then from Lemma~\ref{diam_Santalo}, $D=\diam K$. Suppose for contradiction that $K$ is not complete. Let $\widetilde{K}\supsetneq K$ be a completion of $K$, and $y\in\widetilde{K}\setminus K$. As $K$ is a convex body, there is a supporting hyperplane $H$ to $K$ separating $y$ and $K$. But then, if $z_0\in H\cap\partial K$, since $w_S\left(K,\left(z_0,H\right)\right)=D$, there is a $z'\in\partial K$ and a $H'$ supporting hyperplane to $K$ at $z'$, such that $\left[z_0,z'\right]$ is a diametral chord orthogonal both to $H$ and $H'$. But then
$$
D=\diam\widetilde{K}\geq w_S\left(\widetilde{K},z'\right)>D=w_S\left(K,z'\right)
$$
is a contradiction, therefore $K$ has to be complete. So (iv) implies (i).


We again assume that $K$ is of constant width $D$. Let $H$ be any hyperplane intersecting $K$. Let $A$ and $C$ be the two supporting hyperspheres determined by $H$ (where in the degenerate case one of them is the supporting hyperplane $H$), and let $z\in A\cap\partial K$. If $\ell_z$ is the geodesic line through $z$ and orthogonal to $H$, then the hyperplane $H_z\ni z$ orthogonal to $\ell_z$ supports $K$ at $z$. Let $z'$ be the intersection point of $\ell_z$ and $\partial K\setminus\left\{z\right\}$. If $v\in T_z$ is the outer normal unit vector, such that $-v$ points along the geodesic segment $\left[z,z'\right]$, then $d\left(z,z'\right)=D$ by the assumption. On the other hand, by the orthogonality of $\ell_z$ and $H$, $z'\in C$ follows. Hence, (iii) implies (vii) and (vi).

It is obvious, that (vii) implies (vi). We now assume, that for a fixed interior point $p\in\Int K$ each hyperplane $H$ containing $p$, the Leichtweiss width $w_L\left(K,\left(p,H\right)\right)=D$. We know from Lemma~\ref{diam_Leichtweiss}, that $D=\diam K$. Suppose for contradiction, that $K$ is not complete, so there is a completion $\widetilde{K}\supsetneq K$ of $K$ and a point $y\in\widetilde{K}\setminus K$. But then, there is a boundary point $z\in\partial K$ and a hyperplane $H_z$ which supports $K$ at $z$, separating $K$ and $y$. If $\ell$ is the line through $z$ and orthogonal to $H_z$, and $H$ is the hyperplane containing $p$ and orthogonal to $\ell$, then
$$
D=w_L\left(\widetilde{K},\left(p,H\right)\right)>w_L\left(K,\left(p,H\right)\right)=D
$$
is a contradiction. Therefore, $K$ is already complete, hence (vi) also implies (ii).
\endproof

Additionally, for convex bodies in the hyperbolic plane, Lemma~\ref{wJCJLwS}, Corollary~\ref{diam_JCJL} and Theorem~\ref{completeness} implies the following.

\begin{prop}
Let $K\subset \hyp^2$ a convex body of diameter $D$. 
Then, $K$ is of constant width if and only if $w_{JCJL}\left(K,\left(z,\ell_z,z'\right)\right)=D$ for all $z\in\partial K$, $\ell_z$ supporting line at $z$ and $z'$ opposite point to $z$ with respect to $\ell_z$ 
\end{prop}
\noindent{\bf Remark. }However, being of constant width $D$ is not equivalent with being of constant G. Horv\'ath width in $\hyp^2$. As G. Horv\'ath pointed out in his paper\cite{Hor21}, there are regular $\left(2k+1\right)$-gons of constant G. Horv\'ath width, which are not complete (a completion of theirs are the well-known Reuleaux polygons).

We deduce the following properties of complete sets from Theorem~\ref{completeness} and Lemma~\ref{ballinhoroball}.

\begin{corollary}\label{CompleteProperties}
Let $K\subset \hyp^n$ be a complete set of diameter $D>0$. Then the following hold:
\begin{enumerate}
    \item $K$ is convex and compact;
    \item for all $z\in\partial K$ there is a $y\in\partial K$ such that $d\left(y,z\right)=D$;
    \item for any pair of points $y,z\in K$ if $\sigma$ denotes the shorter arc of some circle of radius at least $D$ through $y$ and $z$, then $\sigma\subseteq K$;
    \item for any pair of points $y,z\in K$ if a horocyclic arc $\theta$ connects $y$ and $z$, then $\theta\subseteq K$.
\end{enumerate}
\end{corollary}

\section{The circumscribed ball and inscribed balls of convex bodies}
\label{balls}

For a compact convex set $K$ in $\hyp^n$, we  define the inradius and circumradius of $K$ as the maximal (resp. minimal) radius of a ball contained in (resp. containing) $K$. We use the notations $r\left(K\right)$ and $R\left(K\right)$, so
$$
r\left(K\right)=\max\left\{r\geq 0\colon\exists z\in \hyp^n\colon B\left(z,r\right)\subseteq K\right\}
$$
and
$$
R\left(K\right)=\min\left\{R\geq 0\colon\exists z\in \hyp^n\colon K\subseteq B\left(z,r\right)\right\}.
$$

We fix a regular $d$-dimensional simplex in $\hyp^n$ of edge length $D$ and denote it by $\Delta\left(D\right)$. Note that $\diam\left(\Delta\left(D\right)\right)=D$.  

The following theorem of Jung is well-known in $\hyp^n$ (see Dekster\cite{Dek95J}), but we need a slightly different version. For the proof of the exact statement below, see B\"or\"oczky, Sagmeister\cite{BoS22}.

\begin{lemma}\label{jung}
If $K\subset \hyp^n$ is a compact set of diameter $D>0$, then
\begin{description}
    \item[(i)] there is a unique point $p$ such that $K\subset B\left(p,R\left(K\right)\right)$;
    \item[(ii)] for some $2\leq k\leq n+1$, there is a set $\left\{q_1,\ldots,q_k\right\}\subset\partial K\cap\partial B\left(p,R\left(K\right)\right)$ such that $p\in\mathrm{relint}\left[q_1,\ldots,q_k\right]$ where $\left[q_1,\ldots,q_k\right]$ is a $\left(k-1\right)$-dimensional simplex in $\hyp^n$;
    \item[(iii)] $$\frac{D}{2}\leq R\left(K\right)\leq R\left(\Delta\left(D\right)\right)<D;$$
    \item[(iv)]  $R\left(K\right)=R\left(\Delta\left(D\right)\right)$ holds if and only if $k=n+1$  and $\left[q_1,\ldots,q_{n+1}\right]$ is congruent with $\Delta\left(D\right)$ in (ii).
\end{description}
\end{lemma}

Concerning the inradius, we have the following result. For a proof, see B\"or\"oczky, Sagmeister\cite{BoS22}.

\begin{lemma}\label{closest-points-on-the-boundary}
If $K\subseteq \hyp^n$ is a convex body, and $B\left(w,r\left(K\right)\right)\subset K$, then for some $2\leq k\leq n+1$ there are points $t_1,\ldots,t_k\in\partial K\cap\partial B\left(w,r\left(K\right)\right)$ such that $\left[t_1,\ldots,t_k\right]$ is a $\left(k-1\right)$-dimensional simplex and $w\in{\rm{relint}}\left[t_1,\ldots,t_k\right]$.
\end{lemma}

\noindent {\bf Remark.} Note that contrary to the circumscribed ball, the inscribed ball of a compact convex set $K$ is not necessarily unique in the Euclidean or the hyperbolic space. For example, take a segment $s$ and let $K$ be the set of points of distance at most $r$ from $s$ for a fixed $r>0$.

However, for complete bodies the following holds (see Böröczky, Sagmeister \cite{BoS22}).

\begin{prop}\label{r_plus_R}
If 
$K\subset \hyp^n$ is a complete set of diameter $D>0$, then
$$
R\left(K\right)+r\left(K\right)=D.
$$
Furthermore, $K$ has a unique inscribed ball whose center is the circumcenter.
\end{prop}

Combining Lemma~\ref{jung} and Proposition~\ref{r_plus_R} implies the following.

\begin{corollary}
\label{CompleteRadii}
If 
$K\subset \hyp^n$ is a complete set of diameter $D>0$, then
$$
R\left(K\right)\leq R\left(\Delta(D)\right) \mbox{ \ and \ }r\left(K\right)\geq D-R\left(\Delta(D)\right).
$$
Furthermore, either $R\left(K\right)= R\left(\Delta(D)\right)$  or $r\left(K\right)= D-R\left(\Delta(D)\right)$
if and only if $K$ contains a congruent copy of
$\Delta(D)$.
\end{corollary}

\noindent{\bf Remark.} According to Zorn's Lemma, there exists a complete set $K_0$ of diameter $D$ in $\hyp^n$ containing a congruent copy of
$\Delta(D)$. If $n=2$, then there exists a unique such set up to congruency; namely, the Reuleaux triangle. For an arbitrary $n\ge2$ the proof of the following theorem gives a construction of such a set.

The following is well-known by Scott \cite{S81} in the Euclidean case, however, the Euclidean proof can be easily extended to $\hyp^n$, and since the original paper is hard to access, we present the proof below.

\vspace{6pt}\noindent{\textbf{Theorem~\ref{scott}.}}
\emph{Let $K\subset \hyp^n$ be a convex body of diameter $D$ and circumscribed ball $B$. Then, there exists a completion $\widetilde{K}$ of $K$ such that $B$ is also the circumscribed ball of $\widetilde{K}$.}
\proof
Let
$$
U\left(K\right)=\left\{x\in B\colon\diam\left(K\cup\left\{x\right\}\right)=D\right\}
$$
and
$$
\varrho\left(K\right)=\sup\left\{d\left(x,K\right)\colon x\in U\left(K\right)\right\}.
$$
Since $U\left(K\right)\subseteq B$, $U\left(K\right)$ is obviously bounded, and it is also easy to see that it is closed: if $\left(y_j\right)_{j\in\N}\subseteq U\left(K\right)$ is a convergent sequence and $y$ is its limit, by the compactness of $B$ we have $y\in B$ and
$$
D\leq\diam\left(K\cup\left\{y\right\}\right)\leq\limsup_{j\in\N}\diam\left(K\cup\left\{y_j\right\}\right)=D,
$$
therefore $y\in U\left(K\right)$. This makes $U\left(K\right)$ compact, and hence
$$
\varrho\left(K\right)=\max\left\{d\left(x,K\right)\colon x\in U\left(K\right)\right\}.
$$
Let $x_1\in U\left(K\right)$ such that $d\left(x_1,K\right)=\varrho\left(K\right)$, and let $K_1=\conv\left(K\cup\left\{x_1\right\}\right)$. We know that $\diam\left(K\cup\left\{x_1\right\}\right)=D$. Furthermore $\diam K_1=D$, since taking the convex hull preserves the diameter (see e.g. B\"or\"oczky, Sagmeister \cite{BoS20} for a reference). Also, $B$ is convex, so $K_1\subseteq B$, and hence we can define $U\left(K_1\right)$ the same way as earlier. Again, this will be a compact set in $B$, and we can define $\varrho\left(K_1\right)$ similarly as before.

For $j\in\N$, $x_{j+1}$ will be a point in $U\left(K_j\right)$, such that $d\left(x_{j+1},K\right)=\varrho\left(K_j\right)$ and let $K_{j+1}=\conv\left(K_j\cup\left\{x_{j+1}\right\}\right)$. Now $\left(K_j\right)_{j\in\N}$ is a bounded sequence of convex bodies in $\hyp^n$ such that $K_j\subseteq K_{j+1}\subseteq B$ for all $j\in\N$, therefore there is a unique convex body $\widetilde{K}\subseteq B$, such that $\widetilde{K}=\lim_{j\to\infty}K_j$. Since $\diam K_j=D$ for all $j\in\N$, we also have $\diam\widetilde{K}=D$. It remains to show that $\widetilde{K}$ is complete.

Let us suppose for contradiction that $\widetilde{K}$ is not complete, and let $y\in \hyp^n\setminus\widetilde{K}$ such that $\diam\left(K\cup\left\{y\right\}\right)=D$. Let $\delta=d\left(y,\widetilde{K}\right)>0$. This also means that $K_{j+1}\setminus K_j\neq\emptyset$, $\varrho\left(K\right),\varrho\left(K_j\right)>0$, and $x_j\neq x_{j+1}$ for all $j\in\N$. We first assume that $y\in B$. Let $1\leq j<k$, hence $K_j\subseteq K_{k-1}$. From the definition of $x_j$ and $\varrho\left(K_{k-1}\right)$, we have that
$$
d\left(x_j,x_k\right)\geq\varrho\left(K_{k-1}\right).
$$
We now show that $\varrho\left(K_{k-1}\right)\geq\delta$. Since $K\subseteq K_{k-1}\subseteq\widetilde{K}$, $y\in B$ and $\diam\left(\widetilde{K}\cup\left\{y\right\}\right)=D$, $\diam\left(K_{k-1}\cup\left\{y\right\}\right)=D$ follows, so $y\in U\left(K_{k-1}\right)$. Also, $K_{k-1}\subseteq\widetilde{K}$ implies $d\left(y,K_{k-1}\right)\geq\delta$. We deduce $\varrho\left(K_{k-1}\right)\geq d\left(y,K_{k-1}\right)\geq\delta$, and hence $d\left(x_j,x_k\right)\geq\delta$, which contradicts the fact that $\left(K_{\alpha}\right)_{\alpha\in\N}$ is convergent. Secondly, let us assume that $y\in \hyp^n\setminus B$. If $z$ is the centre of $B$, $d\left(y,z\right)\geq R+\delta$, where $R>0$ denotes the radius of $B=B\left(z,R\right)$. Let $z_0$ be the unique point on the geodesic line through $y$ and $z$, such that $z\in\left[y,z_0\right]$ and $d\left(y,z_0\right)=D+\delta$ (here we needed $R<D$, which holds by Lemma~\ref{jung}). It is clear from $d\left(y,z\right)=R+\delta$ and $R<D$, that $z_0\not\in\widetilde{K}$, but again from Lemma~\ref{jung} we have that $d\left(z,z_0\right)=D-R\leq R$, so $z_0\in B\setminus\widetilde{K}$. But then $\widetilde{K}\subset B\subset B\left(z_0,D\right)$, so $\diam\left(\widetilde{K}\cup\left\{z_0\right\}\right)=D$ leads to a contradiction, as we have seen it in the previous case. This concludes the proof.
\endproof

\noindent {\bf Remark. }Note that the same proof works in the spherical space with the geodesic metric for complete sets of diameter less than $\frac{\pi}{2}$.

\noindent{\bf Acknowledgement: } 
K\'aroly J. B\"or\"oczky and \'Ad\'am Sagmeister are supported by NKFIH project K 132002.

\end{document}